\documentclass[10pt,twosided]{article}
\usepackage[tbtags]{amsmath}
\usepackage{amssymb}
\allowdisplaybreaks[4]

\usepackage{amssymb,color}

\definecolor{c20}{rgb}{0.,0.7,0.}
\definecolor{c30}{rgb}{0.,0.,1.}
\definecolor{c40}{rgb}{1,0.1,0.7}
\definecolor{c50}{rgb}{1,0,0}
\definecolor{c60}{rgb}{0,0.9,0.1}

\usepackage{amssymb,color}

\newcommand{\abs}[1]{\lvert #1 \rvert}

\newcommand{\E}[1]{\mathbb{E}\left(#1\right)}

\newcommand{\pk}[1]{\mathbb{P} \left( #1 \right) }

\newcommand{\R}{\mathbb{R}}

\newcommand{\limit}[1]{\lim_{#1 \to   \infty}}
\newcommand{\todis}{\stackrel{d}{\to}}

\def\CC{\mathbb{Q}}

\newcommand{\BQN}{\begin{eqnarray}}
\newcommand{\EQN}{\end{eqnarray}}
\newcommand{\BQNY}{\begin{eqnarray*}}
\newcommand{\EQNY}{\end{eqnarray*}}

\newcommand{\BS}{\begin{sat}}
\newcommand{\ES}{\end{sat}}
\newcommand{\BT}{\begin{theo}}
\newcommand{\ET}{\end{theo}}
\newcommand{\BK}{\begin{korr}}
\newcommand{\EK}{\end{korr}}

\newcommand{\BD}{\begin{de}}
\newcommand{\ED}{\end{de}}

\newcommand{\BIT}{\begin{itemize}}
\newcommand{\EIT}{\end{itemize}}
\newcommand{\BDI}{\begin{description}}
\newcommand{\EDI}{\end{description}}

\newcommand{\BRM}{\begin{remarks}}
\newcommand{\ERM}{\end{remarks}}

\newcommand{\BEL}{\begin{lem}}
\newcommand{\EEL}{\end{lem}}

\newtheorem{theo}{Theorem}[section]
\newtheorem{sat}[theo]{Proposition}
\newtheorem{de}[theo]{Definition}
\newtheorem{lem}[theo]{Lemma}

\newtheorem{korr}[theo]{Corollary}

\newtheorem{remarks}[theo]{Remarks}

\newcommand{\nelem}[1]{{Lemma \ref{#1}}}

\newcommand{\netheo}[1]{{Theorem \ref{#1}}}

\newcommand{\prooftheo}[1]{ \textbf{Proof of Theorem} \ref{#1} }

\newcommand{\prooflem}[1]{\textbf{Proof of Lemma} \ref{#1}}

\newcommand{\COM}[1]{}

\newcommand{\QED}{\hfill $\Box$}

\def\Ga{\gamma}

\def\vn{\varepsilon}

\def\rw{\rightarrow}

\def\IF{\infty}
\def\piter{\mathcal{P}}

\date{}

\def\oo{(1+o(1))}

\newcommand{\equaldis}{\stackrel{d}{=}}

\def\CC{\mathbb{Q}}
\def\LT{\left}
\def\RT{\right}
\def\FF{\mathcal{F}}

\def\NN{\mathcal{N}}
\def\ooo{(1+o(1))}

\begin{document}

\title{\bf On the $\gamma$-reflected Processes with fBm Input}

\bigskip
\author{ Peng Liu\thanks{ School Mathematical Sciences and LMPC, Nankai University, Tianjin 300071, China},
Enkelejd Hashorva\thanks{Department of Actuarial Science,
University of Lausanne, UNIL-Dorigny 1015 Lausanne, Switzerland},
 and Lanpeng Ji$^\dagger$}

 \maketitle
\vskip -0.61 cm

  \centerline{\today{}}

\bigskip
{\bf Abstract:} Define a $\gamma$-reflected process $W_\Ga(t)=Y_H(t)-\Ga\inf_{s\in[0,t]}Y_H(s)$, $t\ge0$ with input process $\{Y_H(t), t\ge 0\}$ which is a fractional Brownian motion with Hurst index $H\in (0,1)$ and a negative linear trend. In risk theory
$R_\gamma(t)=u-W_\Ga(t),  t\ge0$ is referred to as the risk process with tax  of a loss-carry-forward type, whereas in queueing theory
$W_1$ is referred to as the queue length process. In this paper, we investigate
the ruin probability and the ruin time  of the risk process $R_\gamma$ over a surplus dependent time interval $[0, T_u]$. 



\section{Introduction}
Let $\{X_H(t), t\ge0\}$ be  a standard fractional Brownian motion (fBm) with Hurst index  $H\in(0,1)$ meaning that $X_H$ is
a centered Gaussian process with almost surely continuous sample paths and covariance  function
\BQNY
Cov(X_H(t),X_H(s))=\frac{1}{2}(\abs{t}^{2H}+\abs{s}^{2H}-\mid t-s\mid^{2H}),\quad t,s\ge0
\EQNY
We define a $\gamma$-reflected process with input process $Y_H(t)= X_H(t)- ct, c>0$ by
\BQN\label{Wgam}
W_\Ga(t)=Y_H(t)-\Ga\inf_{s\in[0,t]}Y_H(s), \  \ t\ge0, \label{PCW}
\EQN
where $\gamma \in [0,1]$ is the reflection parameter. \\
Motivations for studying $W_\Ga$ come from its wide applications in the fields of queuing, insurance, finance and telecommunication. For instance, in queuing theory $W_1$ is referred to as the queue length process (or the workload process); see, e.g., \cite{Harr85, Asm87, ZeeGly00, Whitt02, AwaGly09} among many others.
In  risk  theory the process  $R_\gamma(t)=u-W_\Ga(t),  t\ge0,u\ge 0$ is referred to as
the risk process with tax payments of a loss-carry-forward type; see, e.g., \cite{AsmAlb10}. We refer to \cite{DeMan03, DebickiRol02, HP99, HP08, HZ08, dieker2005extremes} for some recent studies of $W_0$.

For any $u\ge0$, define the {\it ruin time} of the $\gamma$-reflected process $W_\gamma$ by
\BQN\label{eq:tau1}
\tau_{\Ga,u}=\inf\{t\ge0: W_\gamma(t)>u\}\ \ (\text{with}\ \inf\{\emptyset\}=\IF).
\EQN
Further let $T_u, u\ge0$ be a  positive function and define the {\it ruin probability} over a surplus dependent time interval $[0,T_u]$ by
\BQNY
\psi_{\Ga, T_u}(u):=\pk{\tau_{\Ga,u}\le T_u}.
\EQNY
Hereafter,  $\psi_{\Ga, \IF}(u)$ denotes the   ruin probability over an infinite-time horizon.

The ruin time and the ruin probability for the case that $T_u\equiv T\in(0,\IF)$ and the case that $T_u = \IF$ are studied in \cite{HJP13, HJ13}; see also  \cite{DebickiRol02, HP99, HP08, HZ08}. In \cite{HJP13} the exact asymptotics of  $\psi_{\Ga, T}(u)$  and $\psi_{\Ga, \IF}(u)$ are derived, which combined with the results in \cite{HP99} and \cite{DebickiRol02} lead  to the following asymptotic equivalence
\BQN\label{eq:aymequ1}
 \psi_{\Ga, T}(u)=\mathcal{C}_{H, \gamma}\psi_{0,T}(u) \ooo, \ \ \ u\to\IF
 \EQN
 for any $T\in(0,\IF]$, with $\mathcal{C}_{H, \gamma}$ some known positive constant.
The recent contribution  \cite{HJ13} investigates the  approximation of the conditional ruin time $\tau_{\Ga,u} \lvert (\tau_{\Ga,u}<\IF)$.
As shown therein 
the following convergence in distribution (denoted by $\todis$)
  \BQN \label{eq:asmNor1}
\frac{\tau_{\Ga,u}-t_0 u}{A(u)} \Big\lvert (\tau_{\Ga,u}<\IF)\  \todis \ \NN
  \EQN
holds as $u\to\IF$ for any $\Ga\in[0,1)$, where $\NN$ is an $N(0,1)$ random variable and
\BQN\label{t0Au}
t_0=\frac{H}{c(1-H)},\ \ \ A(u)=\frac{H^{H+1/2}u^H}{(1-H)^{H+\frac{1}{2}}c^{H+1}}.
\EQN
See also  \cite{HP08, HZ08, Kobel11, DHJ13a} for related results. We note in passing that the ruin time and the ruin probability are also studied extensively in the framework of other stochastic processes;
see, e.g., \cite{Grif13, GrifMal12,  EKM97, AsmAlb10}.

With motivation from \cite{BoPiPiter09} and \cite{DebKo2013}, as a continuation of the investigation of the aforementioned papers  we shall
analyze the ruin probability and the conditional ruin time of $W_\Ga$ over the surplus dependent time interval $[0,T_u]$ letting $u\to \IF$.
In the literature, commonly the case $T=\IF$ is considered since for many models explicit calculations of the ruin probability is possible.
From a practical point of view, a more interesting quantity is the finite-time ruin probability. The case of the surplus dependent
time interval lies thus in between and is of both practical and theoretical interests; results for the ruin probability in this case (with $\Ga=0$) are initially  derived in \cite{BoPiPiter09}, see Theorem \ref{Coro1} below.
\\
The novel aspect of this paper is that $T_u$ will be a function changing with $u$ according to three different scenarious which we can define with help of (4). 
\COM{following three scenarios (to be specified in Section 2) capturing the behaviour of $T_u$ will be discussed in detail:
\begin{itemize}
\item[i)] The short time horizon;
\item[ii)] The intermediate time horizon;
\item[iii)]The long time horizon.
\end{itemize}

 }
In \netheo{Thm01} below we show that similar asymptotic equivalence as  in \eqref{eq:aymequ1} still holds for all the three scenarios. In \netheo{Thm1} we derive a truncated Gaussian approximation for the (scaled) conditional ruin time over the   long time horizon,
 whereas for the short and the intermediate time horizons 
an exponential approximation is possible. 

Organization of the rest of the paper:  The main results are presented in Section 2 followed then by a section dedicated to the proofs.
In Appendix we present a variant of the celebrated Piterbarg's theorem, which is of interest for further theoretical developments.

\section{Main Results}

In this contribution, three scenarios of $T_u$ will be distinguished. In view of 
\eqref{eq:asmNor1} the asymptotic mean of the ruin time (equal to $t_0 u$)  and the asymptotic standard deviation $A(u)$ (see \eqref{t0Au}) should be used as a scaling parameter for $T_u$ leading to the
definition of the following three scenarios: 
\begin{itemize}
\item[i)] The short time horizon: $\lim_{u\to\IF} T_u/u=0$;
\item[ii)] The intermediate time horizon: $\lim_{u\to\IF} T_u/u=s_0\in(0,t_0)$;
\item[iii)]The long time horizon: $\lim_{u\to\IF}\frac{T_u-t_0u}{A(u)}=x\in (-\IF, \IF]$.
\end{itemize}
Next we introduce two well-known constants appearing in the asymptotic theory of Gaussian processes. Let  therefore $\{B_\alpha(t),t\ge0 \}$ be  a standard fBm  with Hurst index
$\alpha/2 \in (0,1]$.
The {\it Pickands constant} is defined by
\BQNY\label{pick}
\mathcal{H}_\alpha=\lim_{S\rightarrow\infty}\frac{1}{S} \E{ \exp\biggl(\sup_{t\in[0,S]}\Bigl(\sqrt{2}B_\alpha(t)-t^{\alpha}\Bigr)\biggr)}\in(0,\IF),\ \alpha  \in (0,2]
 \EQNY
and the  {\it Piterbarg constant} is  given by
\BQNY
 \mathcal{P}_\alpha^b=\lim_{S\rw\IF}\E{ \exp\biggl(\sup_{t\in[0,S]}\Bigl(\sqrt{2}B_\alpha(t)-(1+b) t ^{\alpha}\Bigr)\biggr)} \in (0,\IF),\ \alpha  \in (0,2],\ b>0.
\EQNY
We refer to \cite{Pit96, debicki2002ruin, DeMan03, debicki2008note, albin2010new, HJP13, DikerY} for  properties and extensions of the Pickands and Piterbarg constants.

In what follows denote by $\Phi(\cdot)$ the distribution function of an $N(0,1)$ random variable, write $\mathcal{N}$ for an $N(0,1)$  random variable
and put $\Psi(\cdot):=1-\Phi(\cdot)$. Before displaying our main results, we include below a key finding of \cite{BoPiPiter09} concerning the ruin probability of the 0-reflected process  $W_0$.

\BT\label{Coro1} Let $W_0$ be the  0-reflected process given as in \eqref{Wgam} with $H\in(0,1)$. We have

i)  If $\lim_{u \rightarrow \infty} {T_u}/{u} = s_0 \in[0, t_0)$, then
\BQN\label{psi01}
\psi_{0, T_u}(u)=D_{H}\left(\frac{u+cT_u}{T^H_u}\right)^{(\frac{1-2H}{H})_{+}}\Psi\left(\frac{u+cT_u}{T^H_u}\right)\ooo
\EQN
 as $u\rw\IF$, where
$$
\mathcal{D}_{H}=\left\{
            \begin{array}{ll}
2^{-\frac{1}{2H}}(H-c_0)^{-1}\mathcal{H}_{{2H}}, & \hbox{if } H <1/2 ,\\
\frac{4(1-c_0)^2}{(1-2c_0)(2-2c_0)}, & \hbox{if } H =1/2,\\
1&  \hbox{if } H >1/2,
              \end{array}
            \right.\ \ \ \text{with}\ c_0=\frac{cs_0}{1+cs_0}.
$$
ii) If $\lim_{u\to\IF}\frac{T_u-t_0u}{A(u)}=x \in (-\IF, \IF] $, then
\BQN\label{psi2}
\psi_{0, T_u}(u)=\psi_{0, \IF}(u)\Phi(x)\ooo
\EQN
  as $u\to\IF$, where the infinite-time ruin probability $\psi_{0,\IF}(u)$ is given by
\BQN
\psi_{0,\IF}(u)=2^{\frac{1}{2}-\frac{1}{2H}}\frac{\sqrt{\pi}}{\sqrt{H(1-H)}} \mathcal{H}_{{2H}}\left(\frac{c^H u^{1-H}}{H^H (1-H)^{1-H}}\right)^{1/H-1}\Psi\left(\frac{c^H u^{1-H}}{H^H (1-H)^{1-H}}\right)(1+o(1)).\label{eq:HP}
\EQN

\ET

The next theorem shows the asymptotic relations between the ruin probability of the  $\Ga$-reflected process $W_\Ga$ and that of the 0-reflected process $W_0$.
Therefore in the light of Theorem \ref{Coro1} we obtain the exact asymptotics of the  ruin probability
over the surplus dependent interval $[0,T_u]$ of the  $\Ga$-reflected process $W_\Ga$.

\BT\label{Thm01} Let $W_\Ga$ be the  $\Ga$-reflected process given as in \eqref{Wgam} with $H\in(0,1)$ and $\Ga\in(0,1)$. We have

 i)  If $\lim_{u \rightarrow \infty} {T_u}/{u} = s_0 \in[0, t_0)$, then
\BQN\label{psi1}
\psi_{\gamma, T_u}(u)=\mathcal{M}_{H,\Ga} \psi_{0, T_u}(u)\ooo
\EQN
 as $u\to\IF$, where
$$
\mathcal{M}_{H,\Ga}=\left\{
            \begin{array}{ll}
\piter_{2H}^{\frac{1-\Ga}{\Ga}}  , & \hbox{if }  H  <1/2 ,\\
\frac{ 2 -2c_0}{2-2c_0-\Ga}, & \hbox{if }  H  =1/2,\\
1&  \hbox{if }  H  >1/2.
              \end{array}
            \right.
$$

ii) If $\lim_{u\to\IF}\frac{T_u-t_0u}{A(u)}=x \in (-\IF, \IF] $, then
\BQN\label{psiGa}
\psi_{\Ga, T_u}(u)=\piter_{2H}^{\frac{1-\Ga}{\Ga}} \psi_{0, T_u}(u)  \ooo
\EQN
 as $u\to\IF$.
\ET

{\bf Remarks.} 
{\it a) For the case that $\Ga=1$ we can add:  Under the statement  $i)$ above
similar arguments as in the proof of \netheo{Thm01} show that \eqref{psi1} holds
 as $u\to\IF$, with
$$
\mathcal{M}_{H,1}=\left\{
            \begin{array}{ll}
 2^{-\frac{1}{2H}}(H-c_0)^{-1}\mathcal{H}_{{2H}}, & \hbox{if }  H  <1/2 ,\\
\frac{ 2 -2c_0}{1-2c_0}, & \hbox{if }  H  =1/2,\\
1&  \hbox{if }  H  >1/2.
              \end{array}
            \right.
$$
For $ii)$, depending on the values of $x$ different asymptotics will appear; those derivations  are more involved and will therefore be omitted here.

b) Another scenario of $T_u$ which is between the cases i) and ii) is that  $\lim_{u\to\IF}\frac{T_u-t_0u}{A(u)}=-\IF.$ This case can not be dealt with in general; 
more conditions should be imposed for the asymptotic behaviour of $T_u$ around $t_0 u$.

c) As discussed in \cite{BoPiPiter09} also of interest is the investigation of the  maximum losses given that ruin occurs,
which, in our setup, is defined as
\BQN\label{defL}
L(\Ga, u):=\left(\sup_{t\in[0,T_u]}W_\Ga(t)-u\right)\Bigg\lvert (\tau_u\le T_u).
\EQN
Under the assumptions of \netheo{Thm01}, we have by an application of \netheo{Coro1} and \netheo{Thm01} that
if $i)$ is satisfied, then
\BQNY
\frac{(u+cT_u) L(\Ga,u)}{T_u^{2H}}  \stackrel{d}{\to}\ \mathcal{E},\ \ \ u\to\IF,
\EQNY
and if $ii)$ is valid, then
 \BQNY
\frac{c^{2H}(1-H)^{2H-1} L(\Ga,u)}{H^{2H}u^{2H-1}}  \stackrel{d}{\to}\ \mathcal{E},\ \ \ u\to\IF.
\EQNY
Here (and in the sequel) $\mathcal{E}$ denotes a unit exponential random variable. Note in passing that the last convergence in distribution is clear when $\Ga=0$ and $T_u=\IF$ since it is known that the random variable $\sup_{t\in[0,\IF)}W_0(t)$ is exponentially distributed with parameter $2c$.}

Below we shall establish asymptotic approximations  for the ruin times considering all three scenarios for $T_u$.
It turns out that for the  long time horizon the (scaled) conditional ruin time can be approximated by a truncated Gaussian random variable.
Surprisingly, this is no longer the case for the short and the intermediate time horizons where the   (scaled) conditional ruin time is approximated by an  exponential random variable.  

\BT\label{Thm1} Let $W_\Ga$ be the  $\Ga$-reflected process given as in \eqref{Wgam} with $H\in(0,1)$ and $\Ga\in(0,1)$,
and let $\tau_{\Ga,u}$ be the ruin time defined as in \eqref{eq:tau1}. We have

i) If $\limit u {T_u}/{u}=0$, then
\BQNY\label{eq:asym1}
\frac{H u^2(T_u-\tau_{\Ga,u})}{T^{2H+1}_u}\Bigl\lvert\LT(\tau_{\Ga,u}\leq T_u\RT)\ \stackrel{d}{\to}\  \mathcal{E},\ \ \ u\to\IF;
\EQNY
ii) If $\limit u {T_u}/{u}=s_0\in(0, t_0)$, then
\BQNY\label{eq:asym2}
\frac{(1+cs_0)(H-(1-H)cs_0)(T_u-\tau_{\Ga,u})}{s_0^{2H+1} u^{2H-1}}\Bigl\lvert\LT(\tau_{\Ga,u} \leq T_u\RT)\ \stackrel{d}{\to}\ \mathcal{E},\ \ \ u\to\IF;
\EQNY
iii) If  $\lim_{u\to\IF}\frac{T_u-t_0u}{A(u)}=x \in (-\IF, \IF]$, then
\BQNY\label{eq:asym3}
\frac{\tau_{\Ga,u}-t_0u}{A(u)}\Bigl\lvert\LT(\tau_{\Ga,u}\leq T_u\RT)\  \stackrel{d}{\to}\ \NN \Bigl\lvert (\NN<x),\ \ \ u\to\IF.
\EQNY

\ET

{\bf Remark.} 
{\it As expected, the  approximation of the conditional ruin time does not involve the reflection constant $\gamma$, since in view of the proof of \netheo{Thm1} the terms with $\Ga$ are canceled out because of the conditional event.
}

\section{Proofs}

In this section, we shall present the proofs of all the theorems. We start with the proof of \netheo{Thm01}. First note that for any $u>0$
\BQNY
\psi_{\gamma, T_u}(u)&=&\mathbb{P}\left(\sup_{t\in[0,T_u]}W_\Ga(t)>u\right)\\
                     &=&\mathbb{P}\left(\sup_{0\leq s \leq t\leq T_u}\Bigl( Z(s,t)-c(t-\Ga s)\Bigr)>u\right),
\EQNY
where $Z(s,t):=X_H(t)-\Ga X_H(s), s,t\ge0.$
Using the self-similarity of the fBm $X_H$, we further have
\BQN\label{eq:psi1}
\psi_{\gamma, T_u}(u)&=&\mathbb{P}\left(\sup_{0\leq s \leq t\leq1}Y_u(s,t)>\frac{u}{T_u^H}\right),
\EQN
where, for any $u>0$
\BQN  \label{eq:Yu}
Y_u(s,t)=\frac{Z(s,t)}{1+\frac{c T_u}{u}(t-\Ga s)},\ \ s,t\ge0.
\EQN

In order to prove statement  $i)$ in \netheo{Thm01}, we give the following crucial lemma.

\BEL\label{LE2}
Let $\{Y_u(s,t),s,t\ge0\}, u>0$ be a family of Gaussian random fields defined as in \eqref{eq:Yu} with $H\in(0,1)$ and $\Ga\in(0,1)$.
Assume that the condition of statement   $i)$ in \netheo{Thm01} is satisfied.
Then, for any $u$ large enough, the variance function $V_{Y_u}^2(s,t) = \E{Y_u^2(s,t)}$ of the Gaussian random field $Y_u$ attains its  maximum over the set $A:=\{(s, t): 0\leq s \leq t \leq 1\}$ at the  unique point $(0,1)$. Moreover,
$$ V_{Y_u}(0,1)=\frac{u}{u+cT_u}.$$
\EEL
\prooflem{LE2} We only present the main ideas of the proof omitting thus some tedious and straightforward calculations.
By solving the two equations
$$\frac{\partial V_{Y_u}^2(s,t)}{\partial s}=0, \quad \frac{\partial V_{Y_u}^2(s,t)}{\partial t}=0$$
 we have that $s=t$. Therefore, the maximum of $V_{Y_u}^2(s,t)$ over $A$ must be attained on the following three lines $l_1=\{(0,t), 0\leq t \leq 1\}$, $l_2=\{(s,t), 0\leq s=t \leq 1\}$ or $l_3=\{(s, 1), 0\leq s \leq 1\}$. It can be shown that on $l_1$ the maximum is attained uniquely at $(0,1)$ and on $l_2$ the maximum is attained uniquely at $(1,1)$. Obviously, both of the two points are on the line $l_3$. Consequently, the  maximum point of $V_{Y_u}^2(s,t)$ over $A$ must be on $l_3$. Moreover, we have that
$$ \left(V_{Y_u}^2(s,1)\right)'= \frac{2c\Ga T_u}{u}\left(1+\frac{cT_u}{u}(1-\Ga s)\right)^{-3}f_{\frac{cT_u}{u}}(s),$$
where, for any $d>0$
\BQN \label{eq:fd}
f_d(s)&=&1-\Ga-(\Ga-\Ga^2)s^{2H}+\Ga(1-s)^{2H}-\frac{H}{d}(1+d-d\Ga s)\nonumber\\
&&\times \left((1-\Ga)s^{2H-1}+(1-s)^{2H-1}\right),\ \ s\ge0.
\EQN
Thus from the following technical lemma   and the fact that
$$\lim_{u\to\IF}\frac{cT_u}{u}=cs_0<\frac{H}{1-H}$$
we conclude that
$$
\left(V_{Y_u}^2(s,1)\right)'<0,\ \ \ \ \forall s\in(0, 1)
$$
implying that the maximum of $V_{Y_u}^2(s,t)$ over the set $A$ is attained at the unique point $(0,1)$. This completes the proof. \QED

\BEL\label{LE1} Let $f_d(s), s\ge0$ be given as in \eqref{eq:fd} with
  $\Ga\in[0, 1)$ and $d\in [0, \frac{H}{1-H})$. Then
  $$
f_d(s)<0,\ \ \ \forall  s\in(0,1).
$$
\EEL
\prooflem{LE1} First rewrite $f_d(s)$  as
\BQNY
f_d(s)&=&(1-\Ga)+\Ga(1-H)(1-s)^{2H}-\Ga(1-\Ga)(1-H)s^{2H}\\
&&-H\LT(1+\frac{1}{d}-\Ga\RT)(1-s)^{2H-1}-\frac{H}{d}(1+d)(1-\Ga)s^{2H-1}.
\EQNY
  Further, we have
  $$1-\Ga<(1-\Ga)(1-s)^{2H-1}+(1-\Ga)s^{2H-1},\ \ s\in(0, 1).$$
   Then, replacing $1-\Ga$ by $(1-\Ga)(1-s)^{2H-1}+(1-\Ga)s^{2H-1}$ in the above equation we obtain
 \BQNY
f_d(s)&<&(1-\Ga)(1-s)^{2H-1}+(1-\Ga)s^{2H-1}+\Ga(1-H)(1-s)^{2H}-\Ga(1-\Ga)(1-H)s^{2H}\\
&&-H\LT(1+\frac{1}{d}-\Ga\RT)(1-s)^{2H-1}-\frac{H}{d}(1+d)(1-\Ga)s^{2H-1}\\
&<&\LT(1-H-\frac{H}{d}\RT)((1-s)^{2H-1}+(1-\Ga)s^{2H-1})-\Ga(1-\Ga)(1-H)s^{2H},
\EQNY
where in the second inequality we used the fact that
$$\Ga(1-H)(1-s)^{2H}\leq\Ga(1-H)(1-s)^{2H-1}, \quad \forall s\in(0,1).$$
Since for any $d\in [0, \frac{H}{1-H})$ 
$$1-H <  \frac{H}{d}$$
we conclude that $f_d(s)<0$ holds for all $s\in(0,1)$, establishing the proof.\QED

\prooftheo{Thm01} $i)$. First, note that (\ref{eq:psi1}) can be rewritten as
$$\psi_{\gamma, T_u}(u)=\mathbb{P}\left(\sup_{0\leq s \leq t\leq1}\frac{Y_u(s,t)}{V_{Y_u}(0,1)}> \frac{u+cT_u}{T_u^H}\right).$$
Next, for any fixed large $u$, we give  expansion of $\frac{V_{Y_u}(s,t)}{V_{Y_u(0,1)}}$ at the point $(0, 1)$. It follows that
\BQN\label{eqcov1}
\frac{V_{Y_u}(s,t)}{V_{Y_u}(0,1)}
&=&\left\{
              \begin{array}{ll}
1-(H-c(u))(1-t)-\Ga(H-c(u))s + o(1-t+s), &  H>1/2 ,\\
1-(\frac{1}{2}-c(u))(1-t)-\Ga (1-\frac{\Ga}{2}-c(u))s+o(1-t+s), & H=1/2,\\
1-(H-c(u))(1-t)-\frac{\Ga-\Ga^2}{2}s^{2H}+o(1-t+s^{2H}), & H<1/2
              \end{array}
            \right.
\EQN
holds as $(s,t)\rw (0,1)$, where $c(u)=\frac{cT_u}{u+cT_u}.$ Furthermore, we have that
\BQN\label{eqcov2}
1-Cov\LT(\frac{Y_u(s,t)}{V_{Y_u}(s,t)}, \frac{Y_u(s',t')}{V_{Y_u}(s',t')}\RT)=\frac{1}{2}\left(\mid t-t'\mid^{2H}+\Ga^2\mid s-s'\mid^{2H}\right)(1+o(1))
\EQN
holds as $(s,t), (s',t')\rw (0,1)$.
In addition, there exists a positive constant $\CC$ such that, for all $u$ large enough
\BQNY
\E{\left(\frac{Y_u(s,t)}{V_{Y_u}(0,1)}-\frac{Y_u(s',t')}{V_{Y_u}(0,1)}\right)^2}\leq
 \CC (|t-t'|^{2H}+|s-s'|^{2H})
\EQNY
holds for all  $(s,t)\in A$. Therefore, by the fact that
$$\lim_{u\to\IF}c(u)=c_0=\frac{cs_0}{1+cs_0}<H$$
 and using \netheo {ThmPiter} (see Appendix), we obtain that
\BQN\label{eq:TuGa}
\psi_{\Ga, T_u}(u)=D_{H,\Ga}\left(\frac{u+cT_u}{T^H_u}\right)^{(\frac{1-2H}{H})_{+}}\Psi\left(\frac{u+cT_u}{T^H_u}\right)\ooo
\EQN
 as $u\rw\IF$, where
$$
\mathcal{D}_{H,\Ga}=\left\{
            \begin{array}{ll}
2^{-\frac{1}{2H}}(H-c_0)^{-1}\mathcal{H}_{{2H}}\piter_{2H}^{\frac{1-\Ga}{\Ga}}, & \hbox{if } H <1/2 ,\\
\frac{4(1-c_0)^2}{(1-2c_0)(2-2c_0-\Ga)}, & \hbox{if } H =1/2,\\
1&  \hbox{if } H >1/2.
              \end{array}
            \right.
$$
Combining the above formula with \eqref{psi01} we obtain (\ref{psi1}). \\
Next, we present the proof of statement $ii)$.\\
 Assume first that $ \lim_{u\to\IF}\frac{T_u-t_0u}{A(u)}=x \in \mathbb{R}$. We have from \eqref{eq:asmNor1} that 
$$\mathbb{P}\LT(\frac{\tau_{\Ga,u}-t_0u}{A(u)}\le x \Bigl\lvert\tau_{\Ga,u}<\IF\RT)\ \rw\ \Phi(x)$$
holds as $u\rw\IF$.
Further note that the above  is equivalent to
$$
\limit{u}\frac{\mathbb{P}\LT(\sup_{0\leq t\leq t_0u+xA(u)}W_\Ga(t)>u\RT)}{\mathbb{P}\LT(\tau_{\Ga,u}<\IF\RT)}= \Phi(x).$$
Thus, we obtain that
$$
\psi_{\Ga,T_u}(u)\ =\ \psi_{\Ga,\IF}(u)\Phi(x) \ooo
$$
as $u\to\IF$, which together with \eqref{psi2} and Theorem 1.1 in \cite{HJP13} yields the validity of  (\ref{psiGa}).
Finally, assume that $ \lim_{u\to\IF}\frac{T_u-t_0u}{A(u)}=\IF$. For any positive large $M$
$$\psi_{\Ga,t_0u+MA(u) }(u)\leq\psi_{\Ga,T_u}(u)\leq\psi_{\Ga,\IF}(u)$$
holds for all $u$ large enough, hence
$$\Phi(M)\leq \liminf_{u\rw\IF}\frac{\psi_{\Ga,T_u}(u)}{\psi_{\Ga,\IF}(u)}\leq\limsup_{u\rw\IF}\frac{\psi_{\Ga,T_u}(u)}{\psi_{\Ga,\IF}(u)}\leq 1.$$
Letting $M\rw\IF$ in the above we conclude that
$$
\psi_{\Ga,T_u}(u)\ =\ \psi_{\Ga,\IF}(u)  \ooo
$$
holds as $u\to\IF$, which further implies that (\ref{psiGa}) is valid, establishing the proof. \QED


 \prooftheo{Thm1} We start with the proof of  statement $i)$. It follows from \eqref{eq:TuGa} that, for any $x>0$
 \BQNY
\mathbb{P}\left(\frac{u^2(T_u-\tau_{\Ga,u})}{T^{2H+1}_u }> x \Bigl\lvert\tau_{\Ga,u}\leq T_u\right)&=&\frac{\mathbb{P}\left(\sup_{0\leq t\leq T_x(u)}W_{\Ga}(t)>u\right)}{\mathbb{P}\left(\sup_{0\leq t\leq T_u}W_{\Ga}(t)>u\right)}\\
&=&\frac{D_{H, \gamma}\left(\frac{u+cT_x(u)}{(T_x(u))^H}\right)^{(\frac{1-2H}{H})_{+}}\Psi\left(\frac{u+cT_x(u)}{(T_x(u))^H}\right)}{D_{H, \gamma}\left(\frac{u+cT_u}{T_u^H}\right)^{(\frac{1-2H}{H})_{+}}\Psi\left(\frac{u+cT_u}{T_u^H}\right)}\oo\\
&=&\exp\LT(-\frac{\LT(\frac{u+cT_x(u)}{(T_x(u))^H}\RT)^2-\LT(\frac{u+cT_u}{T_u^H}\RT)^2}{2}\RT)\oo\\
&\rw&\exp(-H x) 
\EQNY
holds as $u\to\IF$, where $T_x(u)=T_u- {xT^{2H+1}_u }/{u^2}$. Therefore  the claim follows. \\
Next, we give the proof of
statement $ii)$. Similar arguments as above yield  that, for any $x>0$
\BQNY
\mathbb{P}\LT(\frac{T_u-\tau_{\Ga,u}}{u^{2H-1}}>x\Bigl\lvert\tau_{\Ga,u} \leq T_u\RT)&=&\exp\LT(-\frac{\LT(\frac{u+c(T_u-xu^{2H-1})}{(T_u-xu^{2H-1})^H}\RT)^2-\LT(\frac{u+cT_u}{T_u^H}\RT)^2}{2}\RT)\oo\\
&\rw&\exp(-\lambda x) 
\EQNY
holds as $u\to\IF$, where $\lambda=\frac{(1+cs_0)(H-(1-H)cs_0)}{s_0^{2H+1}}$.
Finally, since by \eqref{psiGa} for any $y\le x$
\BQNY
\mathbb{P}\LT(\frac{\tau_{\Ga,u}-t_0u}{A(u)}< y\Bigl\lvert\tau_{\Ga,u}\leq T_u\RT)&=&\frac{\mathbb{P}\LT(\sup_{0\leq t\leq t_0u+yA(u)}W_{\Ga}(t)>u\RT)}{\mathbb{P}\LT(\sup_{0\leq t\leq T_u}W_{\Ga}(t)>u\RT)}\\
& \rw& \frac{\Phi(y)}{\Phi(x)}
\EQNY
holds as $u\to\IF$, the claim of statement $iii)$ follows, and thus the  proof is complete.  \QED

\section{Appendix: Piterbarg's Theorem for Non-homogeneous Gaussian Fields}
We present below a generalization of Theorem  D.3 and Theorem 8.2 in \cite{Pit96}, which is tailored for the proof of the main results. We first introduce a generalization of the
 Piterbarg constant   given by
\BQNY
\tilde{\mathcal{P}}_\alpha^b=\lim_{S\rw\IF}\E{ \exp\biggl(\sup_{t\in[-S,S]}\Bigl(\sqrt{2}B_\alpha(t)-(1+b)\abs{t}^{\alpha}\Bigr)\biggr)} \in (0,\IF),\ \alpha  \in (0,2],\ b>0,
\EQNY
where $\{B_\alpha(t), t\in\R\}$ is a standard fBm defined on $\R$. Let $D=\{(s,t), 0\le s\le t\le 1\}$, and let $\{\eta_u(s,t), (s,t)\in D\}, u\ge0$ be a family of Gaussian random fields satisfying the following three assumptions:

{\bf A1:} The variance function $\sigma_{\eta_u}^2(s,t)$ of $\eta_u$ attains its muximum on the set $D$ at some unique point $(s_0,t_0)$ for any  $u$ large enough, and further there exist four positive constants $A_i, \beta_i, i=1,2$ and two functions $A_i(u), i=1,2$ satisfying $\lim_{u\rw\IF}A_i(u)=A_i, i=1,2$ such that $\sigma_{\eta_u}(s,t)$ has the following expansion around $(s_0,t_0)$ for all $u$ large enough
\BQNY
\sigma_{\eta_u} (s,t)=1-A_1(u)\abs{s-s_0}^{\beta_1}(1+o(1))-A_2(u)\abs{t-t_0}^{\beta_2}(1+o(1)), \ \ (s,t)\to(s_0,t_0).
\EQNY

{\bf A2:} There exist four constants $B_i>0, \alpha_i\in(0,2], i=1,2$ and two functions $B_i(u), i=1,2$ satisfying $\lim_{u\rw\IF}B_i(u)=B_i, i=1,2$ such that the correlation function $r_{\eta_u}(s,t;s',t')$ of $\eta_u$ has the following expansion around $(s_0,t_0)$ for all $u$ large enough
\BQNY
r_{\eta_u} (s,t;s',t')=1-B_1(u)\abs{s-s'}^{\alpha_1}(1+o(1))-B_2(u)\abs{s-s'}^{\alpha_2}(1+o(1)), \ \  (s,t), (s',t')\to (s_0,t_0).
\EQNY

{\bf A3:} For some positive constants $\mathbb{Q}$ and $\gamma$, and all $u$ large enough
$$
\E{\eta_u(s,t)-\eta_u(s',t')}^2\le \mathbb{Q}(\abs{s-s'}^\gamma+\abs{t-t'}^\gamma)
$$
for any $(s,t), (s',t')\in D$.
\BT\label{ThmPiter}
If  $\{\eta_u(s,t), (s,t)\in D\}$, $u\ge0$ is a family of Gaussian random fields satisfying  {\bf A1-A3}, then
$$
\pk{\sup_{(s,t)\in D}\eta_u(s,t)>u}= \FF^{(1)}_{\alpha,\beta}(u)\ \FF^{(2)}_{\alpha,\beta}(u)\ \Psi(u),\ \ \text{as}\ u\rw\IF,
$$
where
$$
\FF^{(i)}_{\alpha,\beta}(u)=\left\{
            \begin{array}{ll}
\widehat{I_i}\mathcal{H}_{\alpha_i} B_i^{\frac{1}{\alpha_i}} A_i^{-\frac{1}{\beta_i}}\Gamma\LT(\frac{1}{\beta_i}+1\RT)\ u^{\frac{2}{\alpha_i}-\frac{2}{\beta_i}}, & \hbox{if } \alpha_i<\beta_i,\\
\widehat{\mathcal{P}}_{\alpha_1}^{\frac{A_i}{B_i}}, & \hbox{if } \alpha_i=\beta_i,\\
1&  \hbox{if } \alpha_i>\beta_i,
              \end{array}
            \right.\ \ \ i=1,2,
$$
with $\Gamma(\cdot)$ the Euler Gamma function and
\BQNY
&&\widehat{\mathcal{P}}_{\alpha_1}^{\frac{A_1}{B_1}}=\left\{
            \begin{array}{ll}
\widetilde{\mathcal{P}}_{\alpha_1}^{\frac{A_1}{B_1}}, & \hbox{if } s_0\in(0,1),\\
\mathcal{P}_{\alpha_1}^{\frac{A_1}{B_1}}&  \hbox{if } s_0=0\ \text{or}\ 1,
              \end{array}
            \right.\ \
 \widehat{\mathcal{P}}_{\alpha_2}^{\frac{A_2}{B_2}}=\left\{
            \begin{array}{ll}
\widetilde{\mathcal{P}}_{\alpha_2}^{\frac{A_2}{B_2}}, & \hbox{if } t_0\in(0,1),\\
\mathcal{P}_{\alpha_1}^{\frac{A_2}{B_2}}&  \hbox{if } t_0=0\ \text{or}\ 1,
              \end{array}
            \right.\\
 &&   \widehat{I_1} =\left\{
            \begin{array}{ll}
2, & \hbox{if } s_0\in(0,1),\\
1&  \hbox{if } s_0=0\ \text{or}\ 1,
              \end{array}
            \right.\ \
 \widehat{I_2} =\left\{
            \begin{array}{ll}
2, & \hbox{if } t_0\in(0,1),\\
1&  \hbox{if } t_0=0\ \text{or}\ 1.
              \end{array}
            \right.
\EQNY

\ET
\prooftheo{ThmPiter} It follows from the assumptions {\bf A1-A2} that for any $\vn>0$ and for $u$ large enough we have
\BQNY
(A_1-\vn)\abs{s-s_0}^{\beta_1}+(A_2-\vn)\abs{t-t_0}^{\beta_2}\le 1-\sigma_{\eta_u} (s,t)\le (A_1+\vn)\abs{s-s_0}^{\beta_1}+(A_2+\vn)\abs{t-t_0}^{\beta_2}
\EQNY
as $ (s,t)\to (s_0,t_0)$, and
\BQNY
(B_1-\vn)\abs{s-s'}^{\alpha_1}+(B_2-\vn)\abs{t-t'}^{\alpha_2}\le 1-r_{\eta_u} (s,t;s',t')\le (B_1+\vn)\abs{s-s'}^{\alpha_1}+(B_2+\vn)\abs{t-t'}^{\alpha_2}
\EQNY
as $ (s,t), (s',t')\to (s_0,t_0)$.
Therefore, in the light of Theorem 8.2 in \cite{Pit96} we can get asymptotical upper and lower bounds, and thus the claims follow by letting $\vn\rw0$. The proof is complete. \QED

\bigskip

{\bf Acknowledgement}: The authors kindly acknowledge partial
support from the Swiss National Science Foundation Project 200021-140633/1,
and the project RARE -318984  (an FP7  Marie Curie IRSES Fellowship).

\bibliographystyle{plain}

 \bibliography{gausbibR}

\end{document}